\documentclass[12pt]{article}
\usepackage{graphicx}
\usepackage{color}
\catcode`\@=11 \@addtoreset{equation}{section}

\catcode`\@=12

\newtheorem{Theorem}{Theorem}[section]
\newtheorem{Proposition}{Proposition}[section]
\newtheorem{Lemma}{Lemma}[section]
\newtheorem{Corollary}{Corollary}[section]
\newtheorem{Remark}{Remark}[section]

\newcommand{\bTheorem}[1]{
\begin{Theorem} \label{T#1} }
\newcommand{\eT}{\end{Theorem}}

\newcommand{\bProposition}[1]{
\begin{Proposition} \label{P#1}}
\newcommand{\eP}{\end{Proposition}}

\newcommand{\bLemma}[1]{
\begin{Lemma} \label{L#1} }
\newcommand{\eL}{\end{Lemma}}

\newcommand{\bCorollary}[1]{
\begin{Corollary} \label{C#1} }
\newcommand{\eC}{\end{Corollary}}

\newcommand{\bRemark}[1]{
\begin{Remark} \label{R#1} }
\newcommand{\eR}{\end{Remark}}

\newcommand{\bFormula}[1]{
\begin{equation} \label{#1}}
\newcommand{\eF}{\end{equation}}

\newcommand{\Ov}[1]{\overline{#1}}

\newcommand{\vr}{\varrho}

\newcommand{\vt}{\vartheta}
\newcommand{\vu}{\vc{u}}
\newcommand{\vc}[1]{{\bf #1}}

\newcommand{\Div}{{\rm div}_x}
\newcommand{\Grad}{\nabla_x}

\newcommand{\dx}{{\rm d} {x}}
\newcommand{\dt}{{\rm d} t }

\newcommand{\intO}[1]{\int_{\Omega} #1 \ \dx}

\definecolor{Cgrey}{rgb}{0.85,0.85,0.85}
\definecolor{Cblue}{rgb}{0.50,0.85,0.85}
\definecolor{Cred}{rgb}{1,0,0}
\definecolor{fancy}{rgb}{0.10,0.85,0.10}

\newcommand\Cbox[2]{%
    \newbox\contentbox%
    \newbox\bkgdbox%
    \setbox\contentbox\hbox to \hsize{%
        \vtop{
            \kern\columnsep
            \hbox to \hsize{%
                \kern\columnsep%
                \advance\hsize by -2\columnsep%
                \setlength{\textwidth}{\hsize}%
                \vbox{
                    \parskip=\baselineskip
                    \parindent=0bp
                    #2
                }%
                \kern\columnsep%
            }%
            \kern\columnsep%
        }%
    }%
    \setbox\bkgdbox\vbox{
        \color{#1}
        \hrule width  \wd\contentbox %
               height \ht\contentbox %
               depth  \dp\contentbox
        \color{black}
    }%
    \wd\bkgdbox=0bp%
    \vbox{\hbox to \hsize{\box\bkgdbox\box\contentbox}}%
    \vskip\baselineskip%
}

\date{}

\begin{document}


\title{Stability of the isentropic Riemann solutions of the full multidimensional Euler system}

\author{Eduard Feireisl \thanks{The research of E.F. leading to these results has received funding from the European Research Council under the European Union's Seventh Framework
Programme (FP7/2007-2013)/ ERC Grant Agreement 320078} \and Ond\v rej Kreml\thanks{O.K. acknowledges the support of the GA\v CR (Czech Science Foundation) project GA13-00522S in the general framework of RVO: 67985840} \and Alexis Vasseur}

\maketitle

\bigskip

\centerline{Institute of Mathematics of the Academy of Sciences of the Czech Republic}

\centerline{\v Zitn\' a 25, 115 67 Praha 1, Czech Republic}

\bigskip

\centerline{Department of Mathematics, University of Texas at Austin}

\centerline{1 University Station C1200, Austin, TX 78712, USA}






\maketitle

\bigskip





\begin{abstract}

We consider the complete Euler system describing the time evolution of an inviscid non-isothermal gas. We show that the rarefaction wave solutions of the 1D Riemann problem are stable, in particular unique, in the class of all bounded weak solutions to the associated multi-D problem. This may be seen as a counterpart of the non-uniqueness results of physically admissible solutions emanating from 1D shock waves constructed recently by the method of convex integration.

\end{abstract}



\section{Introduction}
\label{i}

The recent ground breaking results of De Lellis, Sz\' ekelyhidi and their collaborators \cite{DelSze2}, \cite{DelSze3}, \cite{DelSze} provide a body of evidence that the sofar
well accepted well-posedness criteria for hyperbolic conservation laws based on the Second law of thermodynamics imposed in the form of various entropy
inequalities may fail for certain problems, including the Euler system in gas dynamics, see Chiodaroli et al. \cite{ChiDelKre}. Besides the rather exotic
examples of infinitely many ``wild solutions'' emanating from unspecified and possibly very singular initial data, ill-posedness was demonstrated for
the standard 1-D Riemann data considered in the $N$-D setting, $N=2,3$, see \cite{ChiKre}. Furthermore, some of these solutions dissipate (globally) more kinetic energy than the
standard Riemann solutions; whence the latter apparently violate the entropy rate maximality criterion proposed by Dafermos \cite{Dafer}.

In the light of these examples, it is of interest to study stability of solutions originating from 1-D data considered in the natural multi-D framework. For the mass density $\vr = \vr(t,x)$, the absolute temperature $\vt = \vt(t,x)$, and the velocity field $\vu  = \vu(t,x)$, we introduce the \emph{Euler system}
of partial differential equations:

\Cbox{Cgrey}{

\bFormula{i1}
\partial_t \vr + \Div (\vr \vu) = 0,
\eF
\bFormula{i2}
\partial_t (\vr \vu) + \Div (\vr \vu \otimes \vu) + \Grad (\vr \vt) = 0,
\eF
\bFormula{i3}
\partial_t \left[ \frac{1}{2} \vr |\vu|^2 + c_v \vr \vt \right] + \Div \left[ \left( \frac{1}{2} \vr |\vu|^2 + c_v \vr \vt +
\vr \vt \right) \vu \right] = 0,
\eF
with the associated entropy inequality
\bFormula{i4}
\partial_t (\vr s) + \Div (\vr s \vu) \geq 0, \ s = s(\vr, \vt) \equiv \log \left( \frac{\vt^{c_v}}{\vr} \right).
\eF

}

\noindent
Here $c_v > 0$ denotes the specific heat at constant volume.

For the sake of simplicity, we consider the Cauchy problem for the system (\ref{i1} - \ref{i4}) in the 2-D-case, more specifically in the spatial domain
\[
\Omega = R^1 \times \mathcal{T}^1, \ \mbox{where} \ \mathcal{T}^1 \equiv [0,1] \Big|_{ \{ 0, 1 \} }
\ \mbox{is the ``flat'' sphere,}
\]
meaning all functions of $(t,x_1,x_2)$ considered hereafter are 1-periodic with respect to the second spatial coordinate $x_2$.

We introduce 1-D Riemannian data
\bFormula{i5}
\vr(0,x_1,x_2) = R_0(x_1), \ R_0 = \left\{ \begin{array}{l} R_L \ \mbox{for}\ x_1 \leq 0, \\ \\ R_R \ \mbox{for}\ x_1 > 0, \end{array} \right.
\eF
\bFormula{i6}
\vt(0,x_1,x_2) = \Theta_0(x_1), \ \Theta_0 = \left\{ \begin{array}{l} \Theta_L \ \mbox{for}\ x_1 \leq 0, \\ \\ \Theta_R \ \mbox{for}\ x_1 > 0, \end{array} \right.
\eF
\bFormula{i7}
u^1(0,x_1,x_2) = U_0(x_1), \ U_0 = \left\{ \begin{array}{l} U_L \ \mbox{for}\ x_1 \leq 0, \\ \\ U_R \ \mbox{for}\ x_1 > 0, \end{array} \right. ,
\ u^2(0, x_1,x_2) = 0.
\eF

As is well known, see for instance Chang and Hsiao \cite{ChaHsi}, the Riemann problem (\ref{i1} - \ref{i7}) admits a solution
\[
\vr (t,x) = R(t,x_1) = R(\xi), \ \vt(t,x) = \Theta(t,x) = \Theta(\xi) ,\ \vu(t,x) = [ U(t,x), 0 ] = [U(\xi),0)]
\]
depending solely on the self-similar variable $\xi = \frac{x_1}{t}$. Such a solution is unique in the class of $BV$ solutions of the 1-D problem,
see Chen and Frid \cite{CheFr2}, \cite{CheFr1}.

Our goal is to examine well-posedness of the problem (\ref{i1} - \ref{i7}) in the natural multi-D setting. For the sake of simplicity, we focus on the 2-D case, however, the result remains true for general $N \geq 2$. To see that this is a non-trivial task, we recall that the Riemann problem
for the associated \emph{isentropic} problem \emph{is not} well-posed in the class of admissible weak solutions as soon as the Riemann solution
contains shock(s), see Chiodaroli and Kreml \cite{ChiKre}. Such a non-uniqueness result is obtained in a non-constructive way by means of the method
of convex integration. Solutions produced in this way seem to supply some extra amount of kinetic energy in the system, therefore it is plausible
to expect they may be ruled out as soon as the Riemann solution remains conservative, meaning in the case of rarefaction waves and/or contact discontinuities. Indeed we have shown in \cite{FeKr2014} that the rarefaction wave solution is unique in the class of bounded weak solutions to the
2-D \emph{isentropic} Euler system. In this paper, we extend this result to the complete Euler system (\ref{i1} - \ref{i4}).

Similarly to \cite{FeKr2014}, our approach is based on the application of the relative entropy method proposed by Dafermos \cite{Dafer}. More specifically,
following the ideas of DiPerna \cite{DiP79}, we make use of the associated relative \emph{energy} functional in the spirit of Chen and Frid \cite{CheFr2} that proved to be efficient also in the context of viscous fluids modeled by the  Navier-Stokes-Fourier system \cite{FeiNov10}. In contrast to \cite{CheFr2}, however, the problem must be handled in the Eulerian coordinate system which is the main stumbling block making the
proof definitely more involved also with respect to \cite{FeKr2014}. In particular, there are singular terms sitting on
the wave fronts in the relative energy inequality that must be treated with extra care.

The paper is organized as follows. In Section \ref{m}, we introduce the basic concepts concerning weak solutions to the Euler system, and we state our main result. Section \ref{r} is devoted to the relative energy inequality and its implications on stability of weak solutions. The proof of the main theorem is carried over in Section \ref{p} by means of a careful analysis of the singular terms in the relative energy inequality.

\section{Weak solutions, main result}
\label{m}

We consider general bounded weak solutions to the Euler system satisfying physically relevant non-degeneracy conditions:
\bFormula{m1}
0 <  \vr (t,x) \leq \Ov{\vr},\ 0 <  \vt (t,x) \leq \Ov{\vt},\ |s (\vr, \vt) | < \Ov{s},\
|\vu(t,x)| < \Ov{u} \ \mbox{for a.a.}\ (t,x) \in (0,T) \times \Omega.
\eF
In accordance with the initial conditions
(\ref{i5} - \ref{i7}), we impose the far field conditions:
\bFormula{m2}
\lim_{x_1 \to -\infty} \int_0^T \int_{\mathcal{T}^1} | \vr(t,x_1,x_2)  - R_L | \ {\rm d}x_2 \ \dt = 0 ,\
\lim_{x_1 \to \infty} \int_0^T \int_{\mathcal{T}^1} | \vr(t,x_1,x_2)  - R_R | \ {\rm d}x_2 \ \dt = 0,
\eF
\bFormula{m3}
\lim_{x_1 \to -\infty} \int_0^T \int_{\mathcal{T}^1} | \vt(t,x_1,x_2)  - \Theta_L | \ {\rm d}x_2 \ \dt = 0 ,\
\lim_{x_1 \to \infty} \int_0^T \int_{\mathcal{T}^1} | \vt(t,x_1,x_2)  - \Theta_R | \ {\rm d}x_2 \ \dt = 0,
\eF
\bFormula{m4}
\left\{
\begin{array}{c}
\lim_{x_1 \to -\infty} \int_0^T \int_{\mathcal{T}^1} | u^1(t,x_1,x_2)  - U_L | \ {\rm d}x_2 \ \dt = 0,
\lim_{x_1 \to -\infty} \int_0^T \int_{\mathcal{T}^1} | u^2(t,x_1,x_2) | \ {\rm d}x_2 \ \dt = 0, \\ \\
\lim_{x_1 \to \infty} \int_0^T \int_{\mathcal{T}^1} | u^1(t,x_1,x_2)  - U_R | \ {\rm d}x_2 \ \dt = 0,
\ \lim_{x_1 \to \infty} \int_0^T \int_{\mathcal{T}^1} | u^2(t,x_1,x_2) | \ {\rm d}x_2 \ \dt = 0.
\end{array}
\right\}
\eF

\subsection{Weak formulation}

We say that a trio $[\vr, \vt, \vu]$
is a weak solution of the Euler system (\ref{i1} - \ref{i7}) if the following is satisfied:

\begin{itemize}
\item {\bf Positivity:} $\vr$, $\vt$, and $\vu$ are bounded measurable functions in $[0,T] \times \Omega$ satisfying (\ref{m1}).

\item{\bf Equation of continuity:}

\bFormula{m5}
\intO{ \left[ \vr(\tau,x) \varphi (\tau, x) - \vr_0 (x) \varphi (0, x) \right] }
\eF
\[
= \int_0^\tau \intO{ \Big[ \vr(t,x) \partial_t \varphi(t,x) + \vr \vu (t,x) \cdot \Grad \varphi (t,x) \Big] } \ \dt
\]
for any $0 \leq \tau \leq T$, and any test function $\varphi \in C^1_c([0,T] \times {\Omega})$.

\item {\bf Momentum equation:}

\bFormula{m6}
\intO{ \left[ \vr \vu (\tau, x) \cdot \varphi (\tau, x) - \vr_0 \vu_0 (x) \cdot \varphi (0, x) \right] }
\eF
\[
= \int_0^\tau \intO{ \Big[ \vr \vu (t,x) \cdot \partial_t \varphi (t,x) + \vr [\vu \otimes \vu] (t,x) : \Grad \varphi(t,x) + \vr(t,x)
 \vt (t,x) \Div \varphi
(t,x)\Big] } \ \dt
\]
for any $0 \leq \tau \leq T$, and any $\varphi \in C^1_c([0,T] \times {\Omega};R^2)$.

\item {\bf Total energy equation:}

\bFormula{m7}
\intO{ \left[ \left( \frac{1}{2} \vr |\vu|^2 + c_v\vr \vt \right) (\tau,x) \varphi (\tau, x) - \left( \frac{1}{2} \vr_0 |\vu_0|^2 + c_v\vr_0 \vt_0 \right) (x) \varphi (0, x) \right] }
\eF
\[
= \int_0^\tau \intO{ \Big[ \left( \frac{1}{2} \vr |\vu|^2 + c_v \vr \vt \right) (t,x) \partial_t \varphi(t,x) +
\left( \frac{1}{2} \vr |\vu|^2 + c_v \vr \vt + \vr \vt \right) \vu(t,x) \cdot \Grad \varphi(t,x)  \Big] } \ \dt
\]
for any $0 \leq \tau \leq T$, and any test function $\varphi \in C^1_c([0,T] \times {\Omega})$.

\item {\bf Entropy inequality:}

\bFormula{m8}
\intO{  \left[ \vr b(s(\vr, \vt)) (\tau,x) \varphi (\tau, x) - \vr_0 b(s (\vr_0, \vt_0)) (x) \varphi (0, x) \right] }
\eF
\[
\geq
\int_0^\tau \intO{ \Big[ (\vr b(s(\vr,\vt))
\partial_t \varphi(t,x) + \vr b(s(\vr, \vt)) \vu \cdot \Grad \varphi  \Big] } \ \dt
\]
for any $0 \leq \tau \leq T$, any test function $\varphi \in C^1_c([0,T] \times {\Omega})$, $\varphi \geq 0$, and any $b \in C^1$, $b' \geq 0$.
\end{itemize}

\subsection{Shock-free Riemann solutions to the Euler system}

\label{RS}

We consider the class of the Riemann data producing shock-free solutions, more specifically, solutions that are locally Lipschitz in the \emph{open}
set $(0,T) \times \Omega$. This is the case only if:
\begin{itemize}
\item
the entropy $S$ is \emph{constant} in $[0,T] \times \Omega$;
\item
the density $R$ and the temperature $\Theta$ components of the Riemann solutions are interrelated through
\bFormula{m9}
\Theta = R^{\frac{1}{c_v}} \exp \left( \frac{1}{c_v} S \right);
\eF
\item
the density $R = R(t,x_1)$ and the velocity $U = U(t,x_1)$ represent a rarefaction wave solution of the 1-D \emph{isentropic} system
\bFormula{m10}
\partial_t R + \partial_{{x_1}} (R U ) = 0, \ R \left[ \partial_t U  + U \partial_{x_1} U \right] + \exp\left( \frac{1}{c_v} S \right)
\partial_{x_1} R^{\frac{c_v+1}{c_v}} = 0,
\eF
\end{itemize}
see for instance Chang and Hsiao \cite{ChaHsi}.

\subsection{Uniqueness of the Riemann solution, main result}

We are ready to formulate our main result.

\Cbox{Cgrey}{

\bTheorem{m1}

Let $[\vr,\vt, \vu]$ be a weak solution of the Euler system (\ref{i1} - \ref{i4}) in $(0,T) \times \Omega$ originating from the Riemann data
(\ref{i5} - \ref{i7}) and satisfying the far field conditions (\ref{m2} - \ref{m4}). Suppose in addition that the Riemann data (\ref{i5} - \ref{i7}) give rise to the shock-free solution $[R, \Theta, U]$ of the 1-D Riemann problem specified in Section (\ref{RS}).

Then
\[
\vr = R,\ \vt = \Theta, \ \vu = [U,0] \ \mbox{a.a. in}\ (0,T) \times \Omega.
\]

\eT

}

The remaining part of the paper is essentially devoted to the proof of Theorem \ref{Tm1}.

\section{Relative energy}
\label{r}

The concept of relative entropy/energy goes back to the pioneering work of DiPerna \cite{DiP79} and Dafermos \cite{Dafer}. Here, we adopt the relative energy functional in the Eulerian form introduced in \cite{FeiNov10}:
\bFormula{r1}
\mathcal{E} \left( \vr, \vt, \vu \ \Big| \tilde \vr, \tilde \vt, \tilde \vu \right) =
\intO{ \left[ \frac{1}{2} \vr |\vu - \tilde \vu |^2 + H_{\tilde \vt}(\vr, \vt) - \frac{\partial H_{\tilde \vt}(\tilde \vr, \tilde \vt)}{\partial \vr}
(\vr - \tilde \vr) - H_{\tilde \vt}(\tilde \vr, \tilde \vt) \right] },
\eF
where $H_{\tilde \vt}$ is the ballistic free energy,
\[
H_{\tilde \vt}(\vr, \vt) = \vr \left( c_v \vt - \tilde \vt s(\vr, \vt) \right).
\]
The relative energy functional plays the role of distance between a weak solution of the Euler system and \emph{arbitrary} trio of smooth
test function $[\tilde \vr, \tilde \vt, \tilde \vu]$

The key tool that makes $\mathcal{E}$ extremely useful in the stability problems is the following \emph{relative energy inequality}:
\bFormula{r2}
\left[ \mathcal{E} \left( \vr, \vt, \vu \Big| \tilde \vr, \tilde \vt, \tilde \vu \right) \right]_{t=0}^{t = \tau}
\eF
\[
\leq \int_0^\tau \intO{ \Big[ \vr ( \tilde \vu - \vu) \cdot \partial_t \tilde \vu + \vr
(\tilde \vu - \vu) \otimes \vu : \Grad \tilde \vu + (\tilde \vr \tilde \vt - \vr \vt) \Div \tilde \vu  \Big] } \ \dt
\]
\[
- \int_0^\tau \intO{ \left[ \vr \Big( s(\vr, \vt) - s(\tilde \vr, \tilde \vt) \Big) \partial_t \tilde \vt + \vr
\Big( s(\vr, \vt) - s(\tilde \vr, \tilde \vt) \Big) \vu \cdot \Grad \tilde \vt  \right] } \ \dt
\]
\[
+ \int_0^\tau \intO{ \left[ \left( 1 - \frac{\vr}{\tilde \vr} \right) \partial_t (\tilde \vr \tilde \vt) + \left(\tilde \vu - \frac{\vr}{\tilde \vr}\vu \right) \cdot
\Grad (\tilde \vr \tilde \vt)  \right] }
\ \dt.
\]
The inequality holds for any weak solution $[\vr, \vt, \vu]$ of the Euler system (\ref{i1} - \ref{i7}) satisfying the far field conditions
(\ref{m2} - \ref{m4}) and any trio $[\tilde \vr, \tilde \vt, \tilde \vu]$ of continuously differentiable functions such that
\[
\tilde \vr > 0 , \ \tilde \vt > 0, \ \left\{
\begin{array}{c} \tilde \vr = R_L , \tilde \vt = \Theta_L,\ \tilde u^1 = U_L,\ \tilde u^2 = 0
\ \mbox{whenever}\ x_1 < - A, \\ \\
 \tilde \vr = R_R , \tilde \vt = \Theta_R,\ \tilde u^1 = U_R,\ \tilde u^2 = 0
\ \mbox{whenever}\ x_1 > A
\end{array}
\right\}
\]
for some $A > 0$. The proof of (\ref{r2}) can be done in the same way as in \cite{FeiNov10}, with the straightforward modifications to accommodate the inhomogeneous far-field conditions.

Finally, using a simple density argument, we check without difficulty that the Riemann solution $[R,\Theta, [U,0]]$ can be taken as test functions in
(\ref{r2}) to deduce:
\bFormula{r3}
\left[ \mathcal{E} \left( \vr, \vt, \vu \Big| R, \Theta, [U,0] \right) \right]_{t=0}^{t = \tau}
\eF
\[
\leq \int_0^\tau \intO{ \Big[ \vr ( U - u^1) \partial_t U + \vr
(U - u^1) u^1 \partial_{x_1} U + \left( R \Theta - \vr \vt \right) \partial_{x_1} U  \Big] } \ \dt
\]
\[
- \int_0^\tau \intO{ \left[ \vr \Big( s(\vr, \vt) - S \Big) \partial_t \Theta + \vr
\Big( s(\vr, \vt) - S \Big) u^1 \partial_{x_1} \Theta  \right] } \ \dt
\]
\[
+ \int_0^\tau \intO{ \left[ \left( 1 - \frac{\vr}{R} \right) \partial_t (R \Theta) - \left( U - \frac{\vr}{R} u^1 \right)
\partial_{x_1} (R \Theta)  \right] }
\ \dt.
\]

The proof of Theorem \ref{Tm1} is based on careful analysis of the expression on the right-hand side of (\ref{r3}). Similarly to \cite{FeKr2014},
we show that the expression on the right-hand side of (\ref{r3}) is in fact \emph{non-positive}. The proof, however, is much more involved than in
\cite{FeKr2014} due to the appearance of singular ``cross'' terms.

\section{The proof of Theorem \ref{Tm1}}
\label{p}

We start by observing that
\[
\vr ( U - u^1) \partial_t U + \vr
(U - u^1) u^1 \partial_{x_1} U = - \vr (U - u^1)^2 \partial_{x_1} U + \vr (U - u^1)\left( \partial_t U +  U \partial_{x_1} U \right)
\]
\[
- \vr (U - u^1)^2 \partial_{x_1} U - (U - u^1) \frac{\vr}{R} \partial_{x_1} (R \Theta).
\]
Regrouping terms in (\ref{r3}) we arrive at
\bFormula{p1}
\left[ \mathcal{E} \left( \vr, \vt, \vu \Big| R, \Theta, [U,0] \right) \right]_{t=0}^{t = \tau}
\eF
\[
\leq \int_0^\tau \intO{ \Big[ - \vr (U - u^1)^2 \partial_{x_1} U + (R \Theta - \vr \vt) \partial_{x_1} U  \Big] } \ \dt
\]
\[
- \int_0^\tau \intO{ \left[ \vr \Big( s(\vr, \vt) - S \Big) \left( \partial_t \Theta + U \partial_{x_1} \Theta \right)  + \vr
\Big( s(\vr, \vt) - S \Big) ( u^1 - U ) \partial_{x_1} \Theta  \right] } \ \dt
\]
\[
+ \int_0^\tau \intO{  \left( 1 - \frac{\vr}{R} \right) \Big( \partial_t (R \Theta) + U
\partial_{x_1} (R \Theta) \Big) }
\ \dt.
\]

\subsection{Estimates on the slope of rarefaction waves}

Obviously, the integral on the right-hand side is non-zero only on the non-constant part of the rarefaction wave solution $[R, \Theta, U]$.
Introducing the self-similar variable $\xi = \frac{x}{t}$ we have
\[
R'(\xi) ( U(\xi) - \xi ) = - R(\xi) U'(\xi),
\]
\[
U'(\xi) ( U(\xi) - \xi) + \frac{c_v+1}{c_v} \exp\left(  \frac{1}{c_v} S \right) R^{\frac{1-c_v}{c_v}} R'(\xi) = 0,
\]
from which we readily deduce that
\[
U'(\xi) = 0 \Rightarrow R'(\xi ) = 0 \ \mbox{or, otherwise}\
\left| \frac{ R'(\xi) }{U'(\xi) } \right|^2 = \frac{c_v}{c_v+1} \exp \left( - \frac{1}{c_v} S \right) R^{\frac{2c_v-1}{c_v}}(\xi),
\]
which, combined with
\[
R(\xi) = \exp\left( - S \right) \Theta^{c_v} (\xi),
\]
yields the following conclusion:
Either
\[
\partial_{x_1} U = \partial_{x_1} R = \partial_t R = \partial_{x_1} \Theta = \partial_t \Theta = 0
\]
or
\bFormula{p2}
\left| \frac{ \partial_x \Theta }{ \partial_x U } \right|^2 = \frac{1}{c_v(c_v+1)} \Theta.
\eF
Moreover, as observed in \cite{FeKr2014}, the rarefaction wave solution \emph{always} satisfies
\bFormula{p3}
\partial_{x_1} U \geq 0.
\eF

\subsection{Relative heat}

In accordance with (\ref{p3}), we focus on the case $\partial_{x_1} U > 0$.
Going back to (\ref{p1}), we use (\ref{p2}), (\ref{p3}) to obtain
\[
\vr (s - S)  (u^1 - U) \partial_{x_1} \Theta \leq \vr \frac{1}{4} (s - S)^2 \frac{ | \partial_{x_1} \Theta |^2 }{{ \partial_{x_1} U }} +
\vr | u^1 - U|^2 \partial_{x_1} U.
\]
Thus (\ref{p1}) reduces to
\bFormula{p4}
\left[ \mathcal{E} \left( \vr, \vt, \vu \Big| R, \Theta, [U,0] \right) \right]_{t=0}^{t = \tau}
\eF
\[
\leq \int_0^\tau \intO{  (R \Theta - \vr \vt) \partial_{x_1} U   } \ \dt
\]
\[
- \int_0^\tau \intO{ \left[ \vr \Big( s(\vr, \vt) - S \Big) \left( \partial_t \Theta + U \partial_{x_1} \Theta \right)    \right] } \ \dt
+ \int_0^\tau \intO{  \left( 1 - \frac{\vr}{R} \right) \Big( \partial_t (R \Theta) + U
\partial_{x_1} (R \Theta) \Big) }
\ \dt
\]
\[
+ \int_0^\tau \intO{ \vr \frac{1}{4} (s - S)^2 \frac{ | \partial_{x_1} \Theta |^2 }{{ \partial_{x_1} U }} } \ \dt.
\]

Next, we compute
\bFormula{p5}
(R \Theta  - \vr \vt ) \partial_{x_1} U - \vr \left( s - S \right) \left( \partial_t \Theta + U \cdot \partial_{x_1} \Theta \right) +
\frac{ R - \vr }{R} \left( \partial_t (R \Theta) + U \partial_{x_1} (R \Theta) \right)
\eF
\[
=
(R \Theta  - \vr \vt ) \partial_{x_1} U + \Big[ (R- \vr)  - \vr \left( s - S \right) \Big] \left( \partial_t \Theta + U \cdot \partial_{x_1} \Theta \right) +
\frac{ R - \vr }{R} \Theta \left( \partial_t R  + U \partial_{x_1} R  \right)
\]
\[
=
(R \Theta  - \vr \vt ) \partial_{x_1} U  - \frac{1}{c_v} \Big[ (R- \vr)  - \vr \left( s - S \right) \Big] \Theta \partial_{x_1} U +
(\vr - R) \Theta \partial_{x_1} U
\]
\[
= \vr\left[ \Theta - \vt - \frac{1}{c_v} \left( \frac{R}{\vr} - 1 \right) \Theta + \frac{1}{c_v} (s - S) \Theta     \right] \partial_{x_1} U.
\]

It is worth noting that for  $\Theta = \Theta( V, S)$ expressed as a function of the specific volume $V = \frac{1}{R}$ and the entropy $S$, we get
\[
\Theta(V,S) = \exp \left( \frac{1}{c_v} S \right) V^{- \frac{1}{c_v}} ,\
\partial_V \Theta (V,S) = - \frac{1}{c_v} \Theta R, \ \partial_S \Theta(V,S) = \frac{1}{c_v} \Theta;
\]
whence
\[
\Theta - \vt - \frac{1}{c_v} \left( \frac{R}{\vr} - 1 \right) \Theta + \frac{1}{c_v} (s - S)\Theta
\]
is a \emph{negative-definite} quadratic form in the variables $(\frac{1}{\vr} - \frac{1}{R})$, $s - S$.

\subsection{Non-positivity of the relative energy production}

Comparing (\ref{p4}), (\ref{p5}) and using (\ref{p2}) we conclude that it is enough to show that
the function
\[
F_{R,\Theta} (\vr, s): [\vr, s] \mapsto \left[ \Theta - \vt(\vr,s) - \frac{1}{c_v} \left( \frac{R}{\vr} - 1 \right) \Theta + \frac{1}{c_v} (s - S) \Theta     \right]
+ \frac{\Theta}{4c_v(c_v+1)} (s - S)^2
\]
is non-positive for any choice $\vr > 0$, $s \geq S$.

\bRemark{r1}

Here we point out that $s(\vr, \vt) \geq S$ a.a. in $(0,T) \times \Omega$ as a direct consequence of the entropy inequality (\ref{m8}) and our choice of the initial conditions $s(0, \cdot) = S$.

\eR

Seeing that
\[
\vt(\vr,s) = \exp \left( \frac{1}{c_v} s \right) \vr^{\frac{1}{c_v}}, \ \Theta = \exp \left( \frac{1}{c_v} S \right) R^{\frac{1}{c_v}}
\]
we easily obtain
\[
F_{R, \Theta}(\vr,s) = \Theta \left[ 1 -  \exp \left( \frac{1}{c_v} (s- S) \right) \left( \frac{\vr}{R} \right)^{\frac{1}{c_v}} - \frac{1}{c_v} \left( \frac{R}{\vr} - 1 \right) + \frac{1}{c_v} (s - S) + \frac{1}{4c_v(c_v+1)} (s - S)^2 \right],
\]
and,
introducing the new variables
\[
z = \frac{1}{c_v}(s - S) \geq 0, \ y = \left( \frac{\vr}{R} \right)^{\frac{1}{c_v}} > 0,
\]
we observe that it is enough to examine the function
\[
G(y, z) = \left[ 1 -  \exp \left( z \right) y - \frac{1}{c_v} \left( \frac{1}{y^{c_v}} - 1 \right) + z + \frac{c_v}{4(c_v+1)} z^2 \right].
\]

We start by computing
\[
\partial_y G(y,z) = - \exp(z) + \frac{1}{y^{c_v+1}}, \ \partial_z G(y,z) = - \exp(z) y + 1 + \frac{c_v}{2(c_v+1)}z.
\]

Next, we check easily that
\[
y \mapsto G(y,0) \leq 0 \ \mbox{attaining strong global maximum}\ G(1,0) = 0.
\]

Further observation is that there are no critical points of $G$ in the open set $z > 0$, $y > 0$. Indeed, if
\[
 \exp(z) = \frac{1}{y^{c_v+1}},\   \exp(z) y = 1 + \frac{c_v}{2(c_v+1)}z,
\]
then
\[
\exp \left( \frac{c_v}{c_v+1} z \right) = 1 + \frac{c_v}{2(c_v+1)}z ,
\]
where the last relation holds only if $z = 0$.

The final easy observation is that
\[
G(y,z) \to -\infty \ \mbox{as}\ y \to 0, \ y \to \infty \ \mbox{for any fixed}\ z \geq 0.
\]

Next, suppose that $y \geq 1$. Accordingly, we compute
\[
G(y, z) \leq 1 - y - yz - \frac{1}{2} y z^2  - \frac{1}{c_v} \left( \frac{1}{y^{c_v}} - 1 \right) + z + \frac{c_v}{4(c_v+1)} z^2 \leq 0.
\]

Consequently, it remains to control $G$ for $y \in (0, 1)$ and large $z > 0$. To this end, we fix $z > 0$ and consider the function
\[
y \mapsto G(y,z), \ y \in (0,1].
\]
We already know that $G(1,z) < 0$ and $\lim_{y \to 0} G(y,z) = -\infty$. There is exactly one critical point, namely
\[
y = \exp \left( - \frac{1}{c_v+1} z \right),
\]
whereas the corresponding critical value is
\[
\frac{c_v+1}{c_v} \left( 1 -  \exp \left( \frac{c_v}{c_v+1} z \right) \right) + z + \frac{c_v}{4(c_v+1)} z^2 \leq - z - \frac{c_v}{2(c_v+1)} z^2 + z  + \frac{c_v}{4(c_v+1)} z^2 \leq 0.
\]

Thus we have shown that
\[
F_{R,\Theta} (\vr, s) \leq 0 \ \mbox{for all}\ \vr > 0, \ s \geq S
\]
for any $R, \Theta > 0$. Theorem \ref{Tm1} has been proved.

\end{document}